\documentclass[11pt,leqno,oneside,letterpaper]{amsart}

\usepackage{amssymb,amsmath,amsfonts,latexsym, amscd, amsthm, url} 
\usepackage[letterpaper,left=1.75truein, top=1truein]{geometry}
\usepackage{color}
\usepackage[colorlinks,linkcolor=blue,citecolor=blue]{hyperref}

\usepackage{caption,subcaption,graphicx,epsfig}
\usepackage{tikz-cd,calc}
\usepackage{epsfig}
\usepackage{graphicx}
\usepackage{tikz-cd}
\usepackage{array}
\usepackage{float}
\usepackage{mathtools}
\usetikzlibrary{knots}
\usetikzlibrary{decorations.markings}
\definecolor{light-gray}{gray}{0.92}
\definecolor{ultra-light-gray}{gray}{0.97}

\graphicspath{{images/}}

\newtheorem{theorem}{Theorem}[section]
\newtheorem{lemma}[theorem]{Lemma}
\newtheorem{proposition}[theorem]{Proposition}

\newtheorem{corollary}[theorem]{Corollary} 

\theoremstyle{definition}
\newtheorem{remark}[theorem]{Remark}

\newtheorem{conjecture}[theorem]{Conjecture}

\newtheoremstyle{cases}
  {12pt plus 6 pt}
  {2pt}
  {\bfseries}   
  {}
  {\bfseries}
  {.}
  {.5em}
  {}
\theoremstyle{cases}

\numberwithin{subcase}{case} \numberwithin{subsubcase}{subcase}
\numberwithin{equation}{subsection}

\address{D\'epartement de Math\'ematiques, Universit\'e du Qu\'ebec \`a Montr\'eal, 201 Avenue du Pr\'esident-Kennedy, 
Montr\'eal,
QC H2X 3Y7.} 
\email{ba.idrissa@courrier.uqam.ca}
 
\title{strongly quasipositive quasi-alternating links and Montesinos links}
\author{Idrissa Ba}
\begin{document}

\maketitle

\begin{center}
\today
\end{center}

\begin{abstract} The aim of this article is to give a characterization of strongly quasipositive quasi-alternating 
links and detect new classes of strongly quasipositive Montesinos links and non-strongly
quasipositive Montesinos links. In this direction, we show that, if $L$ is an oriented quasi-alternating link with a
quasi-alternating crossing $c$ such that $L_0$ is alternating (where $L_0$ has the induced orientation), then $L$ is definite if and only 
if it is strongly quasipositive (up to mirroring). We also show that if $L$ is an oriented quasi-alternating link with a quasi-alternating crossing $c$ such that $L_0$ is 
fibred or  more generally has a unique minimal genus Seifert surface (where $L_0$ has the induced orientation), then $L$ is definite 
if and only if it is strongly quasipositive (up to mirroring).
\end{abstract}

\section{introduction}
In this paper we assume, unless otherwise stated, that links are oriented and contained in $\mathbb{S}^3$.
A Seifert surface $S$ is {\it quasipositive} if it is isotopic to a braided Seifert surface $S[b]$ for some positive bandword $b$ (see Rudolph [Ru5]).
A link $L$ is {\it strongly quasipositive} if it has a quasipositive Seifert surface.
Let $L$ be a link. Let $\chi(L)$ ($\chi_4(L)$) be the maximal Euler characteristic number of an oriented 2-manifold without closed
components smoothly embedded in $\mathbb{S}^3$ (in $D^4$ respectively) with boundary $L$.  Let $g(L)$, $g_4(L)$, and $g_4^{top}(L)$ be
the $3$-sphere genus, smooth $4$-ball genus, and topologically (locally flat) $4$-ball genus of $L$, respectively. Let $\sigma$ be the signature of
$L$.

Let $L$ be a link. We say that $L$ is {\it definite} if $g(L)=\frac{1}{2×}(|\sigma(L)| - (m-1))$, where $m$ is the number of components of $L$.

For an oriented quasi-alternating link $L$ with a quasi-alternating crossing $c$, the link $L_0$ is given as in Figure \ref{fig:9}. 
\begin{theorem}\label{lem:5}
 Let $L$ be an oriented quasi-alternating link with a quasi-alternating crossing $c$ such that $L_0$ is alternating.
 Then $L$ is definite if and only if it is strongly quasipositive $($up to mirroring$)$.
\end{theorem}
\begin{theorem}\label{tem:5}
 Let $L$ be an oriented quasi-alternating link with a quasi-alternating crossing $c$ such that $L_0$ is either,
 \begin{enumerate}
 \item fibred, or more generally
 \item has a unique minimal genus Seifert surface.
\end{enumerate}
 Then $L$ is definite if and only if it is strongly quasipositive $($up to mirroring$)$.
\end{theorem}

Based in Theorem \ref{lem:5} and Theorem \ref{tem:5} we made the following conjecture:
\begin{conjecture}\label{conj:1}
Let $L$ be an oriented quasi-alternating link.
 Then $L$ is definite if and only if it is strongly quasipositive $($up to mirroring$)$.
\end{conjecture}
\begin{proposition}\label{propo1}
 If $L$ is an oriented quasi-alternating link which is definite then $L_0$ is also definite.
\end{proposition}

Let $M(e; t_1,\cdots, t_r)$ be a Montesinos knot or link, where $t_i=\frac{\beta_i}{\alpha_i}=[c_1^i,\cdots, c_{m_i}^i]$. We make the following
assumptions below: $\alpha_i>1$ and
 $-\alpha_i< \beta_i<\alpha_i$, with {\rm gcd}$(\alpha_i,\beta_i)=1$. We follow the convention that
\begin{equation*}
[a_1,a_2,\cdots, a_n]=\cfrac{1}{a_1-\cfrac{1}{a_2-\cfrac{1}{a_3-\cdots}}}.
\end{equation*}
 
\begin{proposition}\label{Pthm2}
 Let $L=M(e; t_1,\cdots, t_r)$ be a Montesinos link with $e$ even, $r\geq 3$, $t_i=\frac{\beta_i}{\alpha_i}=[c_1^i,\cdots, c_{m_i}^i]$. 
 \begin{enumerate}
 \item If $e\geq 0$ and $t_i=[2q_1^i, 2s_1^i,2q_2^i, 2s_2^i,\cdots, 2q_{n_i}^i, 2s_{n_i}^i]$ where $q_j^i<0$, $s_j^i< 0$ for any $i=1,\cdots, r$, 
 $j=1,\cdots, n_i$, then $L$ is strongly quasipositive.
 \item If $e\geq 0$ and $t_i=[2q_1^i, 2s_1^i,2q_2^i, 2s_2^i,\cdots, 2q_{n_i}^i]$ where $q_j^i<0$, $s_j^i< 0$ for any $i=1,\cdots, r$, 
 $j=1,\cdots, n_i$ for $q_j^i$, and $j=1,\cdots, n_i-1$ for $s_j^i$ , then $L$ is strongly quasipositive.
 \item If $e\geq 0$ and there exists a unique $t_{i_0}=[2q_1^{i_0}, 2s_1^{i_0},2q_2^{i_0}, 2s_2^{i_0},\cdots, 2q_{n_{i_0}}^{i_0}]$ 
 where $q_j^{i_0}<0$, $s_j^{i_0}< 0$ 
 for any $j=1,\cdots, n_{i_0}$ for $q_j^{i_0}$, and $j=1,\cdots, n_{i_0}-1$ for $s_j^{i_0}$
 and \\ $t_i=[2q_1^i, 2s_1^i,2q_2^i, 2s_2^i,\cdots, 2q_{n_i}^i, 2s_{n_i}^i]$ for any $i\neq i_0$ where $q_j^i<0$, $s_j^i< 0$ for any, 
 $j=1,\cdots, n_i$, then $L$ is strongly quasipositive.
\end{enumerate}

\end{proposition}
\begin{proposition}\label{Pthm3}
 Let $K=M(e; t_1,\cdots, t_r)$ be a Montesinos knot with $e$ even, $r\geq 3$ and $t_i=\frac{\beta_i}{\alpha_i}$.
 Suppose that $\alpha_1$ is even and that, for any $i\geq 2$, $\alpha_i$ is odd and $\beta_i$ is even. Let 
 $t_i=[2c_1^i, 2c_2^i,\cdots, 2c_{m_i}^i]$ be the even continued fraction of $t_i=\frac{\beta_i}{\alpha_i}$, for $i=1,\cdots,r$. Note that
 $m_1$ is odd and $m_i$ is even for $i> 1$.
 If there exists $i_{0}$, $2\leq i_{0}\leq r$ such that $c_1^{i_0}=-c_1^{i_0+1}$, then $K$ is not 
 strongly quasipositive.
\end{proposition}

The paper is organized as follows. In the second section we show Theorem \ref{lem:5}, Theorem \ref{tem:5} and Proposition \ref{propo1}.
In section 3 we prove Proposition \ref{Pthm2} and  Proposition \ref{Pthm3}.

\section*{Acknowledgments}
I would like to thank my PhD supervisor Professor Steven Boyer for drawing my attention
to the topic of the current paper and his consistent encouragement and support. I am also grateful as well to Kenneth Baker for his many
helpful comments and encouragement.

\section{strongly quasipositive quasi-alternating links}
The concept of quasi-alternating links is due to (Ozsv\'ath and Szab\'o, 2005).

The set $\mathcal{Q}$ of {\it quasi-alternating} links is the smallest set of links which satisfies the following properties:
 \begin{enumerate}
 \item the unknot is in $\mathcal{Q}$, 
 \item the set $\mathcal{Q}$ is closed under the following operation. Suppose $L$ is any link which admits a projection with 
 a crossing with the following properties:
 \begin{itemize}
 \item both resolution $L_0$ and $L_\infty$ are in $\mathcal{Q}$, 
 \item det$(L)$ = det$(L_0)$ + det$(L_\infty)$
\end{itemize}

\end{enumerate}
then $L$ is in $\mathcal{Q}$.

Such a crossing $c$ is called quasi-alternating crossing.
\begin{figure}[H]
\centering
\def\svgwidth{0.35\columnwidth}
 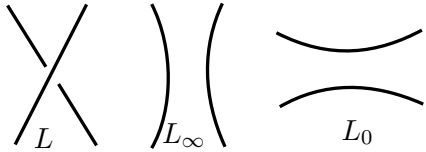
 \caption{The resolutions $L_{\infty}$ and $L_0$ respectively.}
\label{fig:7}
\end{figure}

We have the following concept, which is due to Thurston and taken from [ST].
The {\it complexity} $\chi^-(S)$ of an oriented surface $S$ is $-\chi(C)$, where $C$ is the union of all non-simply connected components of $S$ 
and $\chi(C)$ is the Euler characteristic. For $M$ a compact oriented 3-manifold  and $N$ a (possibly empty) surface in $\partial M$, assign to
any homology class $\alpha\in H_2(M, N,; \mathbb{Z})$ the minimum complexity $x(\alpha)$ of all oriented embedded surfaces whose fundamental
class represents $\alpha$. The function $x: H_2(M, N,; \mathbb{Z})\longrightarrow \mathbb{Z}_+$ is called the {\it Thurston norm}. An oriented
surface ($S, \partial S$)$\subset$($M, \partial M$) is {\it taut} if it is incompressible and 
$\chi^-(S)=x$($[S, \partial S])$ in $H_2(M, \nu(\partial S))$, where $\nu(\partial S)$ is a bicollar neighborhood of $\partial S$ in $\partial M$.

\begin{lemma}{\rm (Lemma in [ST])}\label{lem:2}
 A Seifert surface $S$ for a link $L$ is taut if and only if $\chi(S)=\chi(L)$.
\end{lemma}
\begin{figure}[H]
\centering
\def\svgwidth{0.35\columnwidth}
 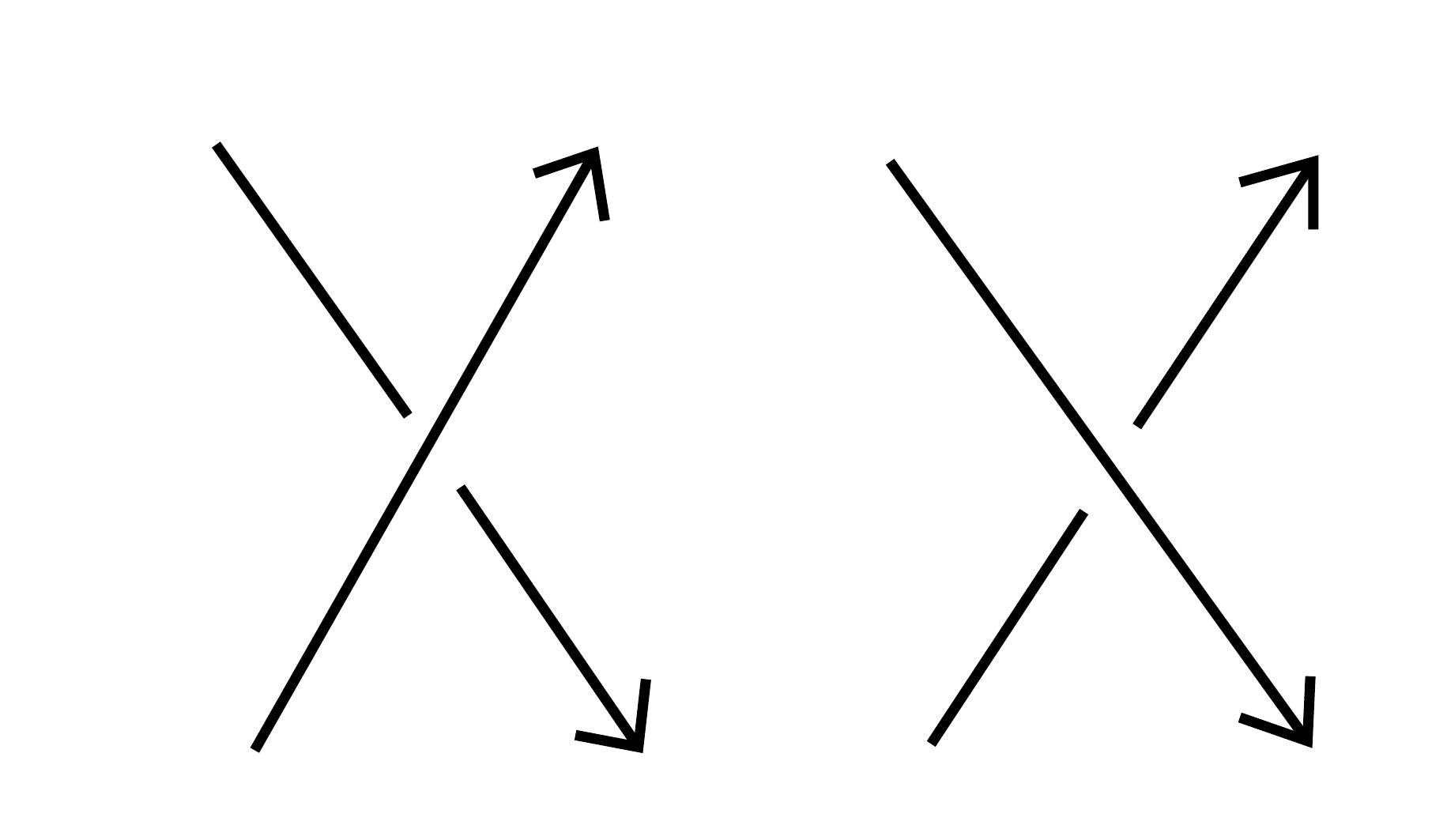
 \caption{A negative crossing and a positive crossing respectively.}
\label{fig:9}
\end{figure}
\begin{corollary}{\rm (Corollary 1.12 in [BBG])}\label{thm: cor2}
 If $L$ is a strongly quasipositive quasi-alternating link then $g_4^{top}(L)=g(L)=\frac{1}{2×}(|\sigma(L)| - (m-1))$, where $m$ is the number 
 of components of $L$.
\end{corollary}

\begin{theorem}{\rm (Lemma 3 in [MO])}\label{thm: them2}
 Let $L$ be an oriented link. The following two conditions are equivalent, providing that the determinants of $L_0$ and $L_{\infty}$ are not zero
 \begin{enumerate}
 \item ${\rm det}(L_{+})={\rm det}(L_{0})+{\rm det}(L_{\infty})$
 \item $\sigma(L_+)=\sigma(L_0)-1$ and $\sigma(L_+)-\sigma(L_{\infty})=-e$.
 
 \end{enumerate}
 A similar equivalence also holds for negative crossings:
\begin{enumerate}
 \item ${\rm det}(L_{-})={\rm det}(L_{0})+{\rm det}(L_{\infty})$
 \item $\sigma(L_-)=\sigma(L_0)+1$ and $\sigma(L_-)-\sigma(L_{\infty})=-e$.
\end{enumerate}
where $e$ is the number of negative crossings in $L_{\infty}$ minus the number of negative crossings in $L_0$. Where also
 $L_+$ and $L_-$ are as in Figure \ref{fig:9}, choose any orientation of $L$ and related orientation of $L_0$ and any orientation of $L_{\infty}$.
\end{theorem}
\begin{figure}[H]
\centering
\def\svgwidth{0.50\columnwidth}
 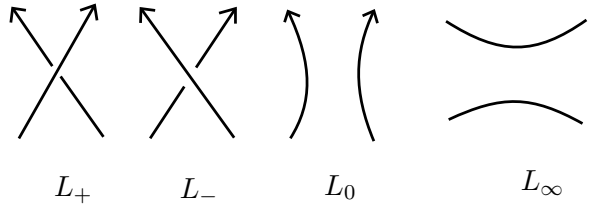
 \caption{The links $L_-$, $L_+$, $L_{\infty}$ and $L_0$ respectively.}
\label{fig:9}
\end{figure}
\begin{theorem}{\rm (Theorem 9.6 in [Li])}\label{thm:17}
Let $L=L_+$ be a link. Then {\rm det} $L_+$ $=$ {\rm det} $L_0$ $+$ {\rm det} $L_{\infty}$ if and only if
 {\rm det} $L_-$ $=$ $|${\rm det} $L_0$ $-$ {\rm det} $L_{\infty}|$, where $L_+$ and $L_-$ are related by a crossing change at a
 crossing.
\end{theorem}
We say that a link is {\it positive} if it admits a diagram with only positive crossings. 
We say that a link is {\it negative} if it admits a diagram with only negative crossings. 
We say that a link is {\it homogeneous} if it has a diagram which can be written as a $\ast$-product of special alternating diagrams (see [Cro]). 
The genus of the Seifert surface obtained by applying the Seifert algorithm to a link diagram $D$ is given by 
$$g(D)=\frac{c(D)-s(D)+2-m}{2},$$ where $c(D)$ is the number of crossings, $s(D)$ is the number of Seifert circles and $m$ is the number of
components of $D$.
\begin{theorem}{\rm (Theorem in [Ru4])}\label{thm:4}
 A Murasugi sum is quasipositive iff the summands are quasipositive.
\end{theorem}
\begin{theorem}{\rm (Theorem 4.1 in [LPS])}\label{thm:3}
 If $D$ is a diagram of a homogeneous link $L$ with $m$ components, then $|w(D)|\leq s(D)+2g(L)+m-1$, where $w(D)$ is the width and $s(D)$ is the
 number of Seifert circle.
 In particular, the equality holds if and only if $D$ is a positive
diagram.
\end{theorem}
\begin{corollary}\label{corl:2}
 If $D$ is a diagram of a homogeneous link $L$, then $g(L)\geq g(D)-c^-(D)$, where $c^-(D)$ is the number of negative crossings in $D$.
\end{corollary}
\begin{theorem}{\rm (Theorem in [Pr1])}\label{thm:100}
 Nontrivial positive links have negative signature.
\end{theorem}
\begin{theorem}{\rm (Theorem 1.2 in [AbT])}\label{them:10}
 A non-split link $L$ is positive if and only if it is homogeneous and strongly quasipositive.
\end{theorem}
\begin{theorem}{\rm (Corollary [Cro])}\label{them:11}
 The surface obtained by applying the Seifert algorithm to a positive diagram of a positive link is a minimal genus Seifert surface.
\end{theorem}

\begin{figure}[H]
\centering
\def\svgwidth{0.50\columnwidth}
 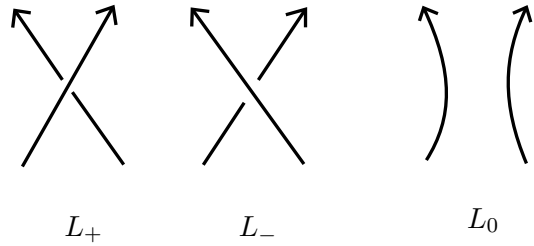
 \caption{We call $L_-$, $L_+$ and $L_0$ Conway moves.}
\label{fig:14}
\end{figure}
\begin{remark}\label{rmk:3}{\rm 1)
 If $L$ is a definite link such that $\sigma(L)=0$ then $L$ is the trivial knot. Indeed, we have that 
 $g(L)=\frac{1}{2×}(|\sigma(L)| - (m-1))=\frac{1}{2×}(0-(m-1))$ which implies that $g(L)=0$ and $m=1$.
 
 2) Figure \ref{fig:19} is obtained by
 considering the Conway moves obtained by giving a half-twist to all the diagrams of Figure \ref{fig:14} (as in Figure \ref{fig:18})
\begin{figure}[H]
\centering
\def\svgwidth{0.60\columnwidth}
 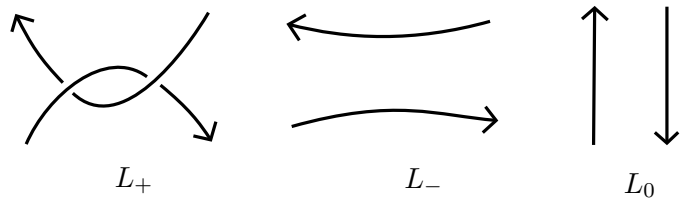
 \caption{ Links obtained by applying half-twist to all the diagrams of Figure \ref{fig:14}.}
\label{fig:18}
\end{figure}
 }
\end{remark}

Let $L$ and $L'$ be links. The distant union of $L$ and $L'$,
written $L\sqcup L'$, is obtained by first moving $L$ and $L'$ so that they are separated by a plane,
and then taking the union.
An arbitrary link $L$ is isotopic to the distant union of its non-splittable sublinks. The number of such non-splittable sublinks is called 
the splitting number.

Since a positive Hopf link is positive then it strongly quasipositive by [Ru2].
\begin{theorem}{\rm (Theorem 1.4 in [ST])}\label{thm:10}
 Let $L_+$, $L_-$ and $L_0$ be three oriented links which are related by the Conway moves at a crossing. Then two of 
 $\chi(L_+)$, $\chi(L_-)$ and $\chi(L_0)-1$ are equal and are no larger than the third. The splitting numbers of the same pair of links are 
 equal and are no larger than the third.
 Therefore, only one of the following three possibilities
 happens,
 \begin{enumerate}
 \item $-\chi(L_+)=-\chi(L_-)\geq 1-\chi(L_0)$;
 \item $-\chi(L_+)=1-\chi(L_0)>-\chi(L_-)$;
 \item $-\chi(L_-)=1-\chi(L_0)>-\chi(L_+)$.
\end{enumerate}
\end{theorem}
\begin{corollary}\label{cor:3}
 Let $L_+$, $L_-$ and $L_0$ be three oriented, non-split links which are related by the Conway moves at a crossing.
 Then only one of the following three possibilities
 happens,
 \begin{enumerate}
 \item $2g(L_+)+ m-1=2g(L_-)+m-1\geq 1+2g(L_0)+m_0-1$;
 \item $2g(L_+)+m-1=1+2g(L_0)+m_0-1>2g(L_-)+m-1$;
 \item $2g(L_-)+m-1=1+2g(L_0)+m_0-1>2g(L_+)+m-1$
\end{enumerate}
where $m$ and $m_0$ are the number of components of $L_{\pm}$ and $L_0$, respectively.
\end{corollary}

\begin{figure}[H]
\centering
\def\svgwidth{0.50\columnwidth}
 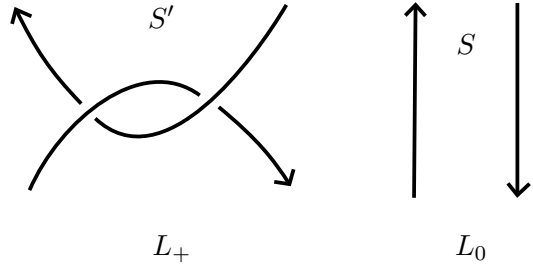
 \caption{ Surfaces $S$ and $S'$.}
\label{fig:10}
\end{figure}

\begin{figure}[H]
\centering
\def\svgwidth{0.60\columnwidth}
 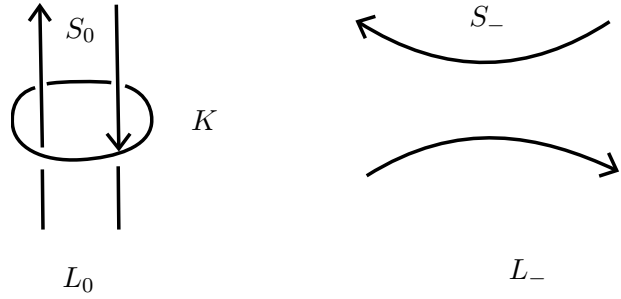
 \caption{ Surfaces $S_0$ and $S_-$.}
\label{fig:19}
\end{figure}

\begin{lemma}{\rm (Proposition 3.1 in [ST])}\label{lem:3}
 When $-\chi(L_+)=1-\chi(L_0)>-\chi(L_-),$ there are taut Seifert surfaces $S'$, $S$ for $L_+$ and $L_0$ respectively
 which appear as in Figure \ref{fig:10}, i.e
 $S'$ is obtained from $S$ by plumbing on a Hopf band: $($An analogous conclusion holds when $-\chi(L_-)=1-\chi(L_0)>-\chi(L_+)$.$)$
\end{lemma}
\begin{proof}[Proof of Proposition \ref{propo1}]Assume that $L$ is a definite, oriented quasi-alternating link.
Since, if $\sigma(L)=0$ then $L$ is the trivial
knot, we can assume that $\sigma(L)\neq0$. Therefore by Theorem \ref{thm: them2} $\sigma(L)$ and $\sigma(L_0)$ will have
the same sign.

 Since $L_+$, $L_-$ and $L_0$ are related by the Conway moves at a crossing, $m_0\in\{m,m-1,m+1\}$. 
 We have that $2g(L_0)\geq |\sigma(L_0)|-(m_0-1)$.
 We have two cases,
 
 \textbf{First case:} If $L=L_+$, then $\sigma(L)=\sigma(L_0)-1$ by Theorem \ref{thm: them2}, and we have that
 \begin{itemize}
 \item $2g(L_0)\geq 2g(L)$ if $\sigma(L)<0$ and $m_0\in\{m,m-1\}$.
 \item $2g(L_0)\geq 2g(L)-2$ if $\sigma(L)<0$ and $m_0=m+1$.
 \item $2g(L_0)> 2g(L)$ if $\sigma(L)>0$ and $m_0\in\{m,m-1\}$.
 \item $2g(L_0)\geq 2g(L)$ if $\sigma(L)>0$ and $m_0=m+1$.
\end{itemize}
Since $L$ and $L_0$ are quasi-alternating links then they are non-split links.
Therefore by (1), (2), (3) of Corollary \ref{cor:3}, Lemma \ref{lem:3} and the fact that $L$ is definite, we have that $m\neq m_0$ and,

$(\ast)$ $\begin{cases} 2g(L_0)=2g(L)$, \hspace{1.1cm} if $m_0=m-1,\; L=L_+\; {\rm and}\; \sigma(L)<0$;
$\\  2g(L_0)=2g(L)-2$,\;\;\; \; if $m_0=m+1 \; L=L_+\; {\rm and}\; \sigma(L)<0$.
$\end{cases}$

In both cases $L_0$ is definite.

 \textbf{Second case:} If $L=L_-$, then $\sigma(L)=\sigma(L_0)+1$ by Theorem \ref{thm: them2}, and we have that
 \begin{itemize}
 \item $2g(L_0)\geq 2g(L)$ if $\sigma(L)>0$ and $m_0\in\{m,m-1\}$.
 \item $2g(L_0)\geq 2g(L)-2$ if $\sigma(L)>0$ and $m_0=m+1$.
 \item $2g(L_0)> 2g(L)$ if $\sigma(L)<0$ and $m_0\in\{m,m-1\}$.
 \item $2g(L_0)\geq 2g(L)$ if $\sigma(L)<0$ and $m_0=m+1$.
\end{itemize}
Since $L$ and $L_0$ are quasi-alternating links then they are non-split links.
Therefore by (1), (2), (3) of Corollary \ref{cor:3}, Lemma \ref{lem:3} and the fact that $L$ is definite, we have that $m\neq m_0$ and, 
 
$(\ast\ast)$ $\begin{cases}  2g(L_0)=2g(L)$, \hspace{1.1cm} if $m_0=m-1,\; L=L_-\; {\rm and}\; \sigma(L)>0$;
$\\  2g(L_0)=2g(L)-2$,\;\;\; \; if $m_0=m+1 \; L=L_-\; {\rm and}\; \sigma(L)>0$. 
$\end{cases}$

In both cases $L_0$ is definite.

Combining $(\ast)$ and $(\ast\ast)$ we have that,

$\begin{cases} 2g(L_0)=2g(L)$, \hspace{1.1cm} if $m_0=m-1,\; L=L_+\; {\rm and}\; \sigma(L)<0$;
$\\ 2g(L_0)=2g(L)$, \hspace{1.1cm} if $m_0=m-1,\; L=L_-\; {\rm and}\; \sigma(L)>0$;
$\\  2g(L_0)=2g(L)-2$,\;\;\; \; if $m_0=m+1 \; L=L_+\; {\rm and}\; \sigma(L)<0$;
$\\  2g(L_0)=2g(L)-2$,\;\;\; \; if $m_0=m+1 \; L=L_-\; {\rm and}\; \sigma(L)>0$. 
$\end{cases}$

\end{proof}
\begin{lemma}{\rm (Claim 2 in [ST])}\label{thm:5}
  Assume that $-\chi(L_+)=1-\chi(L_0)$ and let $S_0$ be any Seifert surface of $L_0$ which is taut in $\mathbb{S}^3\setminus K$ $($as in Figure
  \ref{fig:19}$)$. Then $S_0$ remains taut either in $\mathbb{S}^3$ or in the manifold obtained by doing $0$-surgery along $K$ (see Figure \ref{fig:19}.
  In the latter case, the surface $S_-$ obtained by altering $S_0$ locally as in Figure \ref{fig:19} is a taut Seifert surface 
  for $L_-$ in $\mathbb{S}^3$ and $\chi(S_-)=\chi(S_0)+1$. $($An analogous conclusion holds when $-\chi(L_-)=1-\chi(L_0)$.$)$
\end{lemma}

We have the following theorem which was proved by Boileau-Boyer-Gordon in [BBG] using another method.
\begin{theorem}\label{THEM}
 Let $L$ be a non-split alternating link with $m$ components, then $L$ is definite if and only if it is positive or negative.
\end{theorem}
\begin{proof}
 Let $L$ be an oriented, non-split alternating link. Choose a reduced alternating projection $D$ of $L$. 
Then each crossing of $D$ is a quasi-alternating crossing.

$\Leftarrow$) If the link $L$ is positive or negative then, since $L$ is quasi-alternating, it is definite
by Corollary \ref{thm: cor2} and [Ru2].

$\Rightarrow$) Assume that the link $L$ is definite. We will show by induction on det $L=n$ that $L$ is positive or negative. 
Since, if $\sigma(L)=0$ then $L$ is the trivial
knot, we can assume that $\sigma(L)\neq 0$. Therefore, det $L\geq 2$.

If $n=$ det $L=2$, then, by [LS], $L$ is the Hopf link. Therefore, $L$ is positive or negative.
Now assume that $n= {\rm det}$ $L>2$, and any definite, non-split alternating link $L'$ with det $L'< n$ is positive or negative. We show that 
$L$ is positive or negative.

Since $\sigma(L)\neq 0$, then $\sigma(L)$ and $\sigma(L_0)$ have the same sign by Theorem \ref{thm: them2}. 
Since $L$ is definite then $L_0$ is also definite by
Proposition \ref{propo1}. Therefore, $L_0$ is a positive or negative link by the induction hypothesis or it is the trivial knot.
Without loss of generality assume that $L_0$ is a positive link. Therefore, $\sigma(L_0)\leq 0$ by Theorem \ref{thm:100} and hence $\sigma(L)<0$.
By the proof of Proposition \ref {propo1}, if 
$m$ and $m_0$ are the number of components of $L$ and $L_0$ respectively then
we have that $m\neq m_0$ and

$\begin{cases} 2g(L_0)=2g(L)$, \hspace{1.1cm} if $m_0=m-1,\; L=L_+\; {\rm and}\; \sigma(L)<0$;
$\\ 2g(L_0)=2g(L)$, \hspace{1.1cm} if $m_0=m-1,\; L=L_-\; {\rm and}\; \sigma(L)>0$;
$\\  2g(L_0)=2g(L)-2$,\;\;\; \; if $m_0=m+1 \; L=L_+\; {\rm and}\; \sigma(L)<0$;
$\\  2g(L_0)=2g(L)-2$,\;\;\; \; if $m_0=m+1 \; L=L_-\; {\rm and}\; \sigma(L)>0$.
$\end{cases}$

Therefore, since each crossing of $D$ is a quasi-alternating crossing and $\sigma(L)<0$, each crossing of $D$ is a positive crossing.

\end{proof}
\begin{theorem}{\rm (Proposition 7.1 in [BBG])}
 Let $L$ be a non-split alternating link with $m$ components, then $L$ is definite if and only if it is special alternating.
\end{theorem}
\begin{theorem}\label{them;20}{\rm (Corollary 1.3 in [Gr])} Let $L$ be a special alternating link.
Every minimal genus Seifert surface for $L$ is given by applying Seifert algorithm to some special alternating diagram for $L$.
\end{theorem}

\begin{proof}[Proof of Theorem \ref {lem:5}] Let $L$ be an oriented quasi-alternating link with a quasi-alternating crossing $c$ such that
$L_0$ is alternating.

$\Leftarrow$) If the link $L$ is strongly quasipositive then, since $L$ is quasi-alternating, it is definite
by Corollary \ref{thm: cor2}.

$\Rightarrow$) Assume that the link $L$ is definite. We will show by induction on det $L=n$ that $L$ is strongly quasipositive.
Since, if $\sigma(L)=0$ then $L$ is the trivial
knot, we can assume that $\sigma(L)\neq0$. Therefore, det $L\geq 2$.

If $n=$ det $L=2$, then by [LS] $L$ is the Hopf link. Therefore, $L$ is strongly quasipositive.
Now assume that $n={\rm det}$ $L>2$, and any definite quasi-alternating link $L'$ with a quasi-alternating crossing $c'$ such that
$L_0'$ is alternating with det $L'< n$ is strongly quasipositive. We show that $L$ is strongly quasipositive.

Since $\sigma(L)\neq 0$, then $\sigma(L)$ and $\sigma(L_0)$ have the same sign by Theorem \ref{thm: them2}.
Since $L$ is definite then $L_0$ is also definite by
Proposition \ref {propo1}. Therefore, since $L_0$ is alternating, it is strongly quasipositive by the induction hypothesis.
By the proof of Proposition \ref{propo1}, if $m$ and $m_0$ are the number of components of $L$ and $L_0$, respectively, then
we have that $m\neq m_0$ and

$\begin{cases} 2g(L_0)=2g(L)$, \hspace{1.1cm} if $m_0=m-1,\; L=L_+\; {\rm and}\; \sigma(L)<0$;
$\\ 2g(L_0)=2g(L)$, \hspace{1.1cm} if $m_0=m-1,\; L=L_-\; {\rm and}\; \sigma(L)>0$;
$\\  2g(L_0)=2g(L)-2$,\;\;\; \; if $m_0=m+1 \; L=L_+\; {\rm and}\; \sigma(L)<0$;
$\\  2g(L_0)=2g(L)-2$,\;\;\; \; if $m_0=m+1 \; L=L_-\; {\rm and}\; \sigma(L)>0$.
$\end{cases}$

Therefore, only the following three cases are possible

\begin{enumerate}
  \item $-\chi(L_+)=1-\chi(L_0)>-\chi(L_-)$ with $L=L_+$ and $\sigma(L)<0$;
  \item $-\chi(L_-)=1-\chi(L_0)>-\chi(L_+)$ with $L=L_-$ and $\sigma(L)>0$;
  \item $-\chi(L_+)=-\chi(L_-)= 1-\chi(L_0)$ with $L=L_+$ and $\sigma(L)<0$.
  
\end{enumerate}
We will study each of the cases (1), (2) and (3):

(1) Assume that $-\chi(L_+)=1-\chi(L_0)>-\chi(L_-)$. Then $L=L_+$ and $\sigma(L)<0$.
Therefore by Lemma \ref{lem:3}, there are taut Seifert surfaces $S$, $S_0$ for $L_+$ and $L_0$, respectively,
 which appear as in Figure \ref{fig:10}, i.e $S$ is obtained from $S_0$ by plumbing on a positive Hopf band. Since $L_0$ is definite and
 alternating, it is positive or negative by Theorem \ref{THEM}. Since $\sigma(L)<0$, then $\sigma(L_0)\leq 0$. Therefore $L_0$ is positive.
 Since $L_0$ is positive and alternating then it is special alternating.
 Since $S_0$ is a minimal genus Seifert surface for $L_0$, by Theorem \ref{them;20} $S_0$ is obtained by applying the Seifert algorithm 
 to a special alternating diagram $D_0$ for $L_0$. Since $L_0$ is positive, then $S_0$ is quasipositive by [Ru2].
 Therefore, by Theorem \ref{thm:4} $S$ is also quasipositive, and hence $L$ is strongly quasipositive.
 
(2) Assume that $-\chi(L_-)=1-\chi(L_0)>-\chi(L_+)$. Then $L=L_-$ and $\sigma(L)>0$. In this case, let $L'$ be the mirror image of $L$. Then,
using the same argument as in case (1) for $L'$, we have that $L'$ is strongly quasipositive.
 
(3) Assume that $-\chi(L_+)=-\chi(L_-)= 1-\chi(L_0)$. In this case $L$ can be equal to $L_{+}$ or $L_{-}$. Assume that $L=L_{+}$.
Then $\sigma(L)<0$. By Lemma \ref{thm:5},
if $S_0$ is any Seifert surface of $L_0$ which is taut in $\mathbb{S}^3\setminus K$ (as in Figure
\ref{fig:19}), then we have two subcases:
\begin{itemize}
 \item $S_0$ remains taut in the manifold obtained by doing $0$-surgery to $K$. In this case,
the surface $S_-$ obtained by altering $S_0$ locally as in Figure \ref{fig:19} is a taut Seifert surface for $L_-$ in $\mathbb{S}^3$ and 
$\chi(S_-)=\chi(S_0)+1$. Since $-\chi(L_+)=-\chi(L_-)= 1-\chi(L_0)$ then $S_0$ cannot be taut 
in $\mathbb{S}^3$ because in that case the surface $S_-$ 
obtained by altering $S_0$ locally as in Figure \ref{fig:19} is a taut Seifert surface for $L_-$ in $\mathbb{S}^3$ and 
$\chi(L)=\chi(L_-)=\chi(S_-)=\chi(S_0)+1$. Let $L_2$ be the link obtained
by plumbing on a positive $T(2,4)$ torus link to the Seifert surface $S_0$. Then $|\sigma(L_2)|=|\sigma(L_0)-1|=|\sigma(L)|=2g(L)+m-1$ and 
by Theorem \ref{thm:10} either $\chi(L_2)=\chi(L)=\chi(L_-)$ or $\chi(L_2)=\chi(L_0)-1$ which implies that 
$\chi(L_2)=\chi(L)=\chi(L_-)=\chi(L_0)-1$ and $g(L_2)=g(L)=g(L_-)$, continuing in this way we can plumb on a positive $T(2,2l)$ torus link
and still get a link $L'$ with same genus as $L$. Hence we can assume that $L$ is as in case (1) and then use the same argument.
\item $S_0$ remains taut in $\mathbb{S}^3$.  Since  $-\chi(L)=1-\chi(L_{0})$, the Seifert surface $S$ for $L$ obtained
from $S_0$ by plumbing on a positive Hopf band is taut. Therefore, since $L_0$ is strongly quasipositive and
 alternating, it is special alternating,
 and since $S_0$ is a minimal genus Seifert surface for $L_0$, by Theorem \ref{them;20} $S_0$ is obtained by applying the Seifert algorithm 
 to a special alternating diagram $D_0$ for $L_0$. Since $L_0$ is definite, then $S_0$ is quasipositive by [BBG] and [Ru2].
 Therefore, by Theorem \ref{thm:4} $S$ is also quasipositive, and hence $L$ is strongly quasipositive.
\end{itemize}
The case where $L=L_{-}$ is done in the same way using (2).

This completes the proof that $L$ is strongly quasipositive.
\end{proof}
\begin{theorem}{\rm (Theorem 5.1 in [Ko])}\label{tem:3}
Let $L_1$, $L_2$ be links with minimal genus Seifert surfaces $R_1$, $R_2$, respectively, and $R$ a
Murasugi sum of $R_1$ and $R_2$. Then the minimal genus Seifert surfaces for $L = \partial R$ are
unique if and only if one of $L_1$, $L_2$, say $L_1$, is fibred and the minimal genus Seifert
surfaces for $L_2$ are unique.
\end{theorem}

\begin{proof}[Proof of Theorem \ref{tem:5}] Since fibred links have unique minimal genus Seifert surfaces, it suffices to show
this Theorem when $L_0$ has a unique minimal genus Seifert surface.

Let $L$ be an oriented quasi-alternating link with a quasi-alternating crossing $c$ such that
$L_0$ has a unique minimal genus Seifert surface.

$\Leftarrow$) If the link $L$ is strongly quasipositive, then it is definite by Corollary \ref{thm: cor2}.

$\Rightarrow$) Assume that the link $L$ is definite. We will show by induction on det $L=n$ that $L$ is strongly quasipositive. 
We can assume that $\sigma(L)\neq 0$. Therefore, det $L\geq 2$.

If $n=$ det $L=2$, then by [LS] $L$ is the Hopf link. Therefore, $L$ is strongly quasipositive.
Now assume that $n={\rm det}$ $L>2$, and any definite quasi-alternating link $L'$ with a quasi-alternating crossing $c'$ such that $L_0'$ has a unique
minimal genus Seifert surface and with det $L'< n$ is strongly quasipositive. We show that $L$ is strongly quasipositive. 

Since $L$ is definite then $L_0$ is also definite by
Proposition \ref{propo1}. Hence, if $\sigma(L_0)=0$ then $L_0$ is the trivial knot. Therefore, we can assume that $\sigma(L_0)\neq 0$. Let $c_0$ be a
quasi-alternating crossing for $L_0$. By the proof of Proposition \ref{propo1}, if $m_0$ and $m_{00}$ are the number of components of 
$L_0$ and $L_{00}$, respectively, then we have that $m_0\neq m_{00}$ and

$\begin{cases} 2g(L_{00})=2g(L_0)$, \hspace{1.1cm} if $m_{00}=m_0-1,\; L_0=L_{0+}\; {\rm and}\; \sigma(L_0)<0$;
$\\ 2g(L_{00})=2g(L_0)$, \hspace{1.1cm} if $m_{00}=m_0-1,\; L_0=L_{0-}\; {\rm and}\; \sigma(L_0)>0$;
$\\  2g(L_{00})=2g(L_0)-2$,\;\;\; \; if $m_{00}=m_0+1 \; L_0=L_{0+}\; {\rm and}\; \sigma(L_0)<0$;
$\\  2g(L_{00})=2g(L_0)-2$,\;\;\; \; if $m_{00}=m_0+1 \; L_0=L_{0-}\; {\rm and}\; \sigma(L_0)>0$.
$\end{cases}$

Therefore, only the following three cases are possible

\begin{enumerate}
  \item $-\chi(L_{0+})=1-\chi(L_{00})>-\chi(L_{0-})$ if $L_0=L_{0+}$ and $\sigma(L_0)<0$;
  \item $-\chi(L_{0-})=1-\chi(L_{00})>-\chi(L_{0+})$ if $L_0=L_{0-}$ and $\sigma(L_0)>0$;
  \item $-\chi(L_{0+})=-\chi(L_{0-})= 1-\chi(L_{00})$ if $L_0=L_{0+}$ and $\sigma(L_0)<0$ or $L_0=L_{0-}$ and $\sigma(L_0)>0$.
\end{enumerate}
Hence, we have that:

(1) Assume that $-\chi(L_{0+})=1-\chi(L_{00})>-\chi(L_{0-})$. Then $L_0=L_{0+}$ and $\sigma(L_0)<0$.
Therefore, by Lemma \ref{lem:3}, there are taut Seifert surfaces $S$, $S_0$ for $L_{0+}$ and $L_{00}$, respectively,
which appear as in Figure \ref{fig:10}, i.e $S$ is obtained from $S_0$ by plumbing on a positive Hopf band. Since $L_0$ has a unique minimal
genus Seifert surface, by Theorem \ref{tem:3} $L_{00}$ has also a unique minimal genus Seifert surface.
 
(2) Assume that $-\chi(L_{0-})=1-\chi(L_{00})>-\chi(L_{0+})$. Then $L_0=L_{0-}$ and $\sigma(L_0)>0$. Therefore
by Lemma \ref{lem:3}, there are taut Seifert surfaces $S$, $S_0$ for $L_0$ and $L_{00}$, respectively,
which appear as in Figure \ref{fig:10}, i.e $S$ is obtained from $S_0$ by plumbing on a negative Hopf band. Since $L_0$ has a unique minimal
genus Seifert surface, by Theorem \ref{tem:3} $L_{00}$ has also a unique minimal genus Seifert surface.

(3) Assume that $-\chi(L_{0+})=-\chi(L_{0-})= 1-\chi(L_{00})$. In this case $L_0$ can be equal to $L_{0+}$ or $L_{0-}$. Assume that $L_0=L_{0+}$.
Then $\sigma(L_0)<0$. By Lemma \ref{thm:5}, if $S_0$ is any Seifert surface of $L_{00}$ which is taut in $\mathbb{S}^3\setminus K$ (as in Figure
\ref{fig:19}), then we have two subcases:
\begin{itemize}
 \item $S_0$ remains taut in the manifold obtained by doing $0$-surgery to $K$. In this case,
the surface $S_-$ obtained by altering $S_0$ locally as in Figure \ref{fig:19} is a taut Seifert surface for $L_{0-}$ in $\mathbb{S}^3$ and 
$\chi(S_-)=\chi(S_0)+1$. Since $-\chi(L_{0+})=-\chi(L_{0-})= 1-\chi(L_{00})$ then $S_0$ cannot be taut 
in $\mathbb{S}^3$ because in that case the surface $S_-$ 
obtained by altering $S_0$ locally as in Figure \ref{fig:19} is a taut Seifert surface for $L_{0-}$ in $\mathbb{S}^3$ and 
$\chi(L_0)=\chi(L_{0-})=\chi(S_-)=\chi(S_0)+1$. Let $L_2$ be the link obtained
by plumbing a positive $T(2,4)$ torus link to the Seifert surface $S_0$. Then $|\sigma(L_2)|=|\sigma(L_{00})-1|=|\sigma(L_0)|=2g(L_0)+m_0-1$ and 
by Theorem \ref{thm:10} either $\chi(L_2)=\chi(L_0)=\chi(L_{0-})$ or $\chi(L_2)=\chi(L_{00})-1$ which implies that 
$\chi(L_2)=\chi(L_0)=\chi(L_{0-})=\chi(L_{00})-1$ and $g(L_2)=g(L_0)=g(L_{0-})$, continuing in this way we can plumb on a positive $T(2,2l)$ 
torus link and still get a link $L'$ with same genus as $L_0$. Hence we can assume that $L_0$ is as in case (1) and then use the same argument.
\item $S_0$ remains taut in $\mathbb{S}^3$. Since  $-\chi(L_{0+})=1-\chi(L_{00})$, the Seifert surface $S$ for $L_0$ obtained
from $S_0$ by plumbing on a positive Hopf band is taut. Since $L_0$ has a unique minimal
genus Seifert surface, by Theorem \ref{tem:3} $L_{00}$ has also a unique minimal genus Seifert surface.
\end{itemize}
The case where $L_0=L_{0-}$ is done in the same way using (2).

Therefore, the link $L_0$ has a quasi-alternating crossing $c_0$ such that $L_{00}$ has a unique minimal genus Seifert surface. Thus
$L_0$ is a strongly quasipositive link by the induction hypothesis.
Since $L_0$ has a unique minimal genus Seifert surface then, using a similar argument as in the 
proof of Theorem \ref{lem:5}, we show that $L$ is also strongly quasipositive.
\end{proof}
\begin{remark}\label{rmk:2}{\rm
Since every quasi-alternating link gives rise to an infinite family of quasi-alternating links obtained by replacing
a quasi-alternating crossing with an alternating rational tangle [ChK], 
Theorem \ref{lem:5} and Theorem \ref{tem:5} can be extended to more classes of
quasi-alternating
links.}
\end{remark}
Let L be a link and $D$ a link diagram of $L$. Checkerboard color the regions of the 
complement of the diagram in $\mathbb{R}^2$. Assume that the unbounded region 
$X_0$ is
colored white. The other white regions will be called by $X_1$, $X_2$, $\cdots$, $X_n$. To any crossing $p$ of L we 
associate the number $\chi(p)$ which is +1 or -1 according to the convention in the Figure 1.
\begin{figure}[H]
\centering
\def\svgwidth{0.35\columnwidth}
 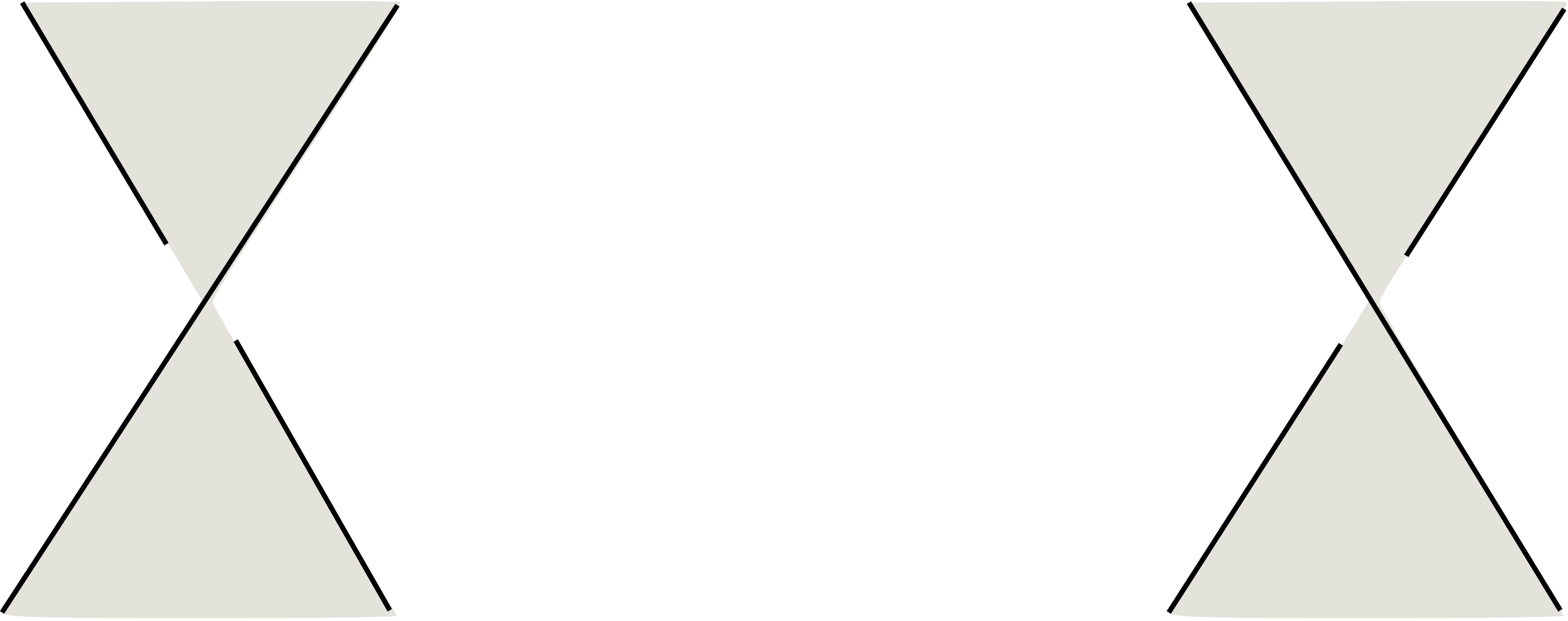
 \caption{The sign convention +1 and -1 respectively.}
\label{figure;6}
\end{figure}

Let $B(D)$ be the following graph: a planar graph whose vertices are in one-to-one correspondence with 
the black regions in our checkerboard coloring, and whose edges correspond to crossings in the diagram. Call this graph 
the black graph of $D$.

\begin{theorem}{\rm (4.6 in [Cr])}\label{thm: main resultP2}
Let $Tr(D)$ be the set of maximal trees in $B(D)$ {\rm(that is trees which contain all the vertices of $B(D)$)}
and $E$ be the set of all edges of $B(D)$. Then 
$${\rm det}\; L=|\sum_{T\in Tr(D)}\prod_{e\in T\cap E}\chi(e)|$$
\end{theorem}
\begin{proposition}\label{thm: main resultP1}
Let $L$ be a non-alternating quasi-alternating link diagram with a quasi-alternating crossing c such 
that the links $L^c_0$ and $L^c_{\infty}$ are alternating. If the Conway notation of $L$ has only one negative sign, then $c$ is positive.
\end{proposition}
\begin{proof}

Assume that the Conway notation of 
$L$ has only one negative sign. By contradiction let $c$ be the negative crossing and $A$ be the set of all maximal trees of $B(D)$ which
contains $c$. Then by Theorem \ref{thm: main resultP2},
we have also that,
\begin{equation} \label{eq1}
\begin{split}
 {\rm det}\; L &= |\sum_{T\in Tr(D)}\prod_{e\in T\cap E}\chi(e)|\\
 &= |\sum_{T\in A}\prod_{e\in A\cap E}\chi(e)+\sum_{T\in Tr(D)\setminus A}\prod_{e\in T\cap E}\chi(e)|\\
 &= |-card(A)+card(Tr(D)\setminus A)|\\
 &<|card(A)+card(Tr(D)\setminus A)|\\
 &< {\rm det}\; L^c_{\infty}+{\rm det}\;L^c_0={\rm det}\;L
\end{split}
\end{equation}
Therefore, the crossing $c$ is positive.
\end{proof}
\begin{remark}
Proposition \ref{thm: main resultP1} states that for any non-alternating, quasi-alternating diagram 
 $D$, the diagrams $D_0$ and $D_\infty$ are not
alternating. They can be alternating only up to isotopy.
\end{remark}
\begin{proposition}\label{thm: main prop2}
Let $L$ be non-alternating quasi-alternating link diagram with a quasi-alternating crossing c such 
that $L^c_0$ and $L^c_{\infty}$ are alternating (up to isotopy) and
whose checkerboard color has only 
two negative crossings of the form {\rm [-2]} or 
{\rm [-$\frac{1}{2}$]}.
If $L=L_{{\rm [-\frac{1}{2}]}}$ and $c$ is negative then $L_{{\rm [-2]}}$ is not quasi-alternating at the two negative crossings
and vice versa, where $L_{{\rm [-2]}}$ is the
link diagram obtained from  $L_{{\rm [-\frac{1}{2}]}}$ by replacing the elementary tangle {\rm [-$\frac{1}{2}$]} by {\rm [-2]}.
\end{proposition}
\begin{proof}
Assume that $L=L_{{\rm [-\frac{1}{2}]}}$ and $c$ is negative. We show that $L'=L_{{\rm [-2]}}$ is not quasi-alternating at 
the two negative crossings. 
Let $c_1$ and $c_2$ be the two negative crossings of $L_{{\rm [-2]}}$.
Without loss of generality we can assume that the path $c_1c_2$ form a loop in $B(L')$.
Any maximal tree of $B(L')$ contains at most one of them.
Let $A$ be the set of all maximal trees of $B(L')$ which contains $c_1$ and $A'$ be the set of all maximal trees of $B(L')$ which contains $c_2$.
We have that
\begin{equation} \label{eq1}
\begin{split}
 {\rm det}\; L' &= |\sum_{T\in Tr(L')}\prod_{e\in T\cap E}\chi(e)|\\
 &= |\sum_{T\in A}\prod_{e\in A\cap E}\chi(e)+\sum_{T\in A'}\prod_{e\in A'\cap E}\chi(e)+
 \sum_{T\in Tr(L')\setminus \{A\cup A'\}}\prod_{e\in T\cap E}\chi(e)|\\
 &= |-card(A)-card(A')+card(Tr(L')\setminus \{A\cup A'\})|\\
 &= |-2{\rm det}\; L'_{\infty}+{\rm det}\;L'_{0,0}|
\end{split}
\end{equation}
Here the resolutions are applied in the two negative crossings in both $L'=L_{{\rm [-2]}}$ and $L=L_{{\rm [-\frac{1}{2}]}}$.

In a similar way, we show that ${\rm det}\; L= |-2{\rm det}\;L_{0}+{\rm det}\; L_{\infty,\infty}|$, ${\rm det}\; L'_0= {\rm det}\; L_{\infty}
= |-{\rm det}\; L_{\infty,\infty}+{\rm det}\;L_{0}|$, ${\rm det}\; L'_{\infty}={\rm det}\; L_{\infty,\infty}$ and 
${\rm det}\; L_0={\rm det}\;L'_{0,0}$.
Since $L=L_{{\rm [-\frac{1}{2}]}}$ is quasi-alternating and $c$ is negative then 
${\rm det}\;L={\rm det}\; L_{\infty}+{\rm det}\;L_0=|-{\rm det}\; L_{\infty,\infty}+{\rm det}\;L_{0}| + {\rm det}\;L_0=
|-2{\rm det}\;L_{0}+{\rm det}\; L_{\infty,\infty}|$, so
${\rm det}\;L=2{\rm det}\;L_{0}-{\rm det}\; L_{\infty,\infty}$, which implies that ${\rm det}\; L_{\infty} = 
{\rm det}\;L_{0} -{\rm det}\; L_{\infty,\infty}$.
Hence, if $L'=L_{{\rm [-2]}}$ was quasi-alternating, then 
${\rm det}\;L'={\rm det}\; L'_{\infty}+{\rm det}\;L'_0={\rm det}\; L_{\infty,\infty} + {\rm det}\;L_{0} -
{\rm det}\; L_{\infty,\infty}= {\rm det}\;L_{0}={\rm det}\;L'_{0,0}$, 
and this contradict the fact that ${\rm det}\; L'=|-2{\rm det}\; L'_{\infty}+{\rm det}\;L'_{0,0}|$.

\end{proof}

\section{strongly quasipositive montesinos links}

\begin{figure}[H]
\centering
\def\svgwidth{0.60\columnwidth}
 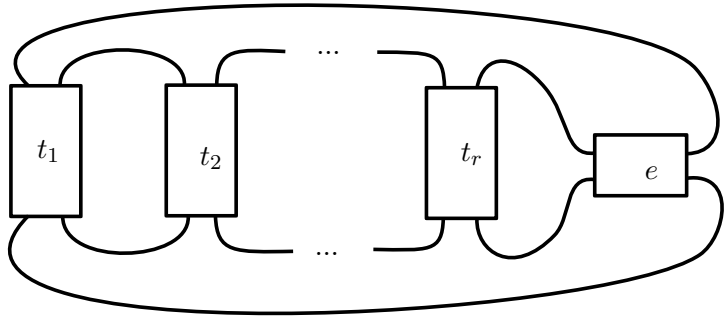
 \caption{A Montesinos link.}
\label{fig:2}
\end{figure}

\begin{figure}[H]
\centering
\def\svgwidth{0.60\columnwidth}
 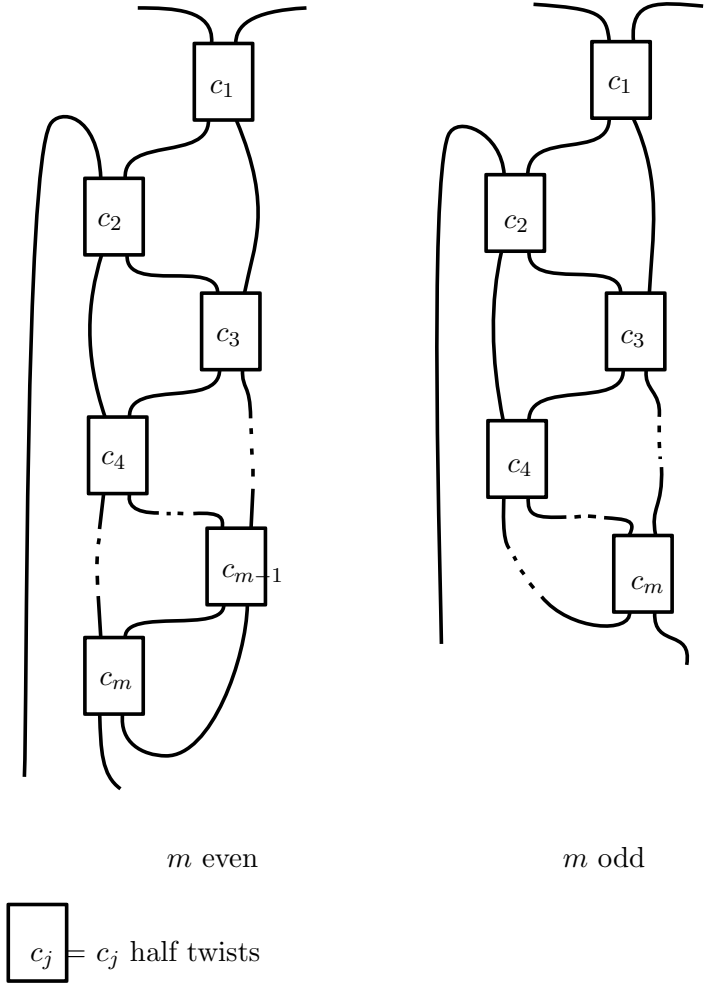
 \caption{The tangle $T(\frac{\beta}{\alpha})$.}
\label{fig:3}
\end{figure}

A Montesinos link, denoted $M(e; t_1,t_2,\cdots, t_r)$, is a link having a diagram of the form depicted in Figure \ref{fig:2},
where $t_i=\frac{\beta_i}{\alpha_i}=[c_1^i,\cdots, c_{m_i}^i]$, each $t_i$ represents a rational tangle (for the definition of rational tangle see
Figure \ref{fig:3}), and $e\in\mathbb{Z}$ represents $|e|$-half twists.
We make the following assumptions in all this paper: $\alpha_i>1$ and $-\alpha_i< \beta_i<\alpha_i$, with {\rm gcd}$(\alpha_i,\beta_i)=1$.\\
If $T(\frac{\beta}{\alpha})$ is a rational tangle of slope $\frac{\beta}{\alpha}$ with
$\alpha>1$, $-\alpha< \beta<\alpha$ and {\rm gcd}$(\alpha,\beta)=1$, then the upper left and the lower left corners of $T(\frac{\beta}{\alpha})$
are connected by a strand if and only if $\alpha$ is even. Therefore, for $K=M(e; t_1,t_2,\cdots, t_r)$ to be a knot rather than a link, at 
most one of $\alpha_1$, $\alpha_2$,$\cdots$, $\alpha_r$ can be even. Since a cyclic permutation of indices does not change the knot type, we assume
that $\alpha_2$,$\cdots$, $\alpha_r$ are odd if $K=M(e; t_1,t_2,\cdots, t_r)$ is a knot.

Let $T(\frac{\beta}{\alpha})$ be a rational tangle of slope $\frac{\beta}{\alpha}$ with
$\alpha>1$, $-\alpha< \beta<\alpha$ and {\rm gcd}$(\alpha,\beta)=1$, such that
\begin{equation*}
 \frac{\beta}{\alpha}=[c_1,c_2,\cdots, c_m]=\cfrac{1}{c_1-\cfrac{1}{c_2-\cfrac{1}{c_3-\cdots}}}.
\end{equation*}

The tangle $T(\frac{\beta}{\alpha})$ is represented in Figure \ref{fig:3}.

Therefore, up to now on each of the $t_i$ box of the Montesinos link $M(e; t_1,t_2,\cdots, t_r)$ in Figure \ref{fig:2}
will be represented by a rational tangle as the one in Figure \ref{fig:3}.
\begin{remark}{\rm 
If $r\leq 2$ then the Montesinos link $M(e; t_1,t_2,\cdots, t_r)$ is a 2-bridge link. Therefore, throughout this paper we can assume that 
$r\geq 3$.
 }
\end{remark}
Let $C=[c_1, c_2,\cdots, c_m]$ be a continued fraction. We say that $C$ is an {\it even continued fraction} if all the $c_j$ are even. We say that
$C$ is a strict continued fraction if
\begin{enumerate}
 \item $c_j$ is even for any odd $j$, and
 \item $c_jc_{j+1}< 0$ whenever $j$ is odd and $|c_j|=2$.
\end{enumerate}
\begin{lemma}{\rm (Proposition 2.5 in [HM])}\label{thm: lem1}
 Let $\frac{\beta}{\alpha}$ be a rational number such that $\alpha$ is odd, $-\alpha< \beta<\alpha$, {\rm gcd}$(\alpha,\beta)=1$, and
 $2|\beta|<\alpha$. Then $\frac{\beta}{\alpha}$ has a strict continued fraction.
\end{lemma}
\begin{figure}[H]
\centering
\def\svgwidth{0.35\columnwidth}
 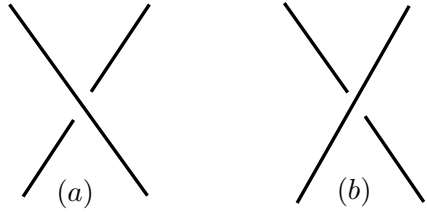
 \caption{Type of crossings.}
\label{fig:4}
\end{figure}
\begin{proof}[Proof of Proposition \ref{Pthm2}] 
(1) In this case $L$ is a 2-component link. Since for any $i=1,\cdots, r$,
$t_i=[2q_1^i, 2s_1^i,2q_2^i, 2s_2^i,\cdots, 2q_{n_i}^i, 2s_{n_i}^i]$ where $q_j^i<0$, $s_j^i< 0$, then each crossing of each
tangle $T(t_i=\frac{\beta_i}{\alpha_i})$ has the form (a) as in Figure \ref{fig:4}. Since $e\geq 0$ then each crossing of $e$ has 
the form (b) as in Figure \ref{fig:4}. Since $e$ is even there exists an orientation of $L$ which makes it positive. Therefore, by [Ru2] $L$ is 
strongly quasipositive.

(2) In this case $L$ is a $r$-components link. Since for any $i=1,\cdots, r$,
$t_i=[2q_1^i, 2s_1^i,2q_2^i, 2s_2^i,\cdots, 2q_{n_i}^i]$ where $q_j^i<0$, $s_j^i< 0$, then each crossing of each
tangle $T(t_i=\frac{\beta_i}{\alpha_i})$ has the form (a) as in Figure \ref{fig:4}. Since $e\geq 0$ then each crossing of $e$ has 
the form (b) as in Figure \ref{fig:4}. Since $e$ is even, there exists an orientation of $L$ which makes it positive. Therefore, by [Ru2] $L$ is 
strongly quasipositive.

(3) In this case $L$ is a knot. Since all the continued fractions are even and $q_j^i<0$, $s_j^i< 0$, then each crossing of each
tangle $T(t_i=\frac{\beta_i}{\alpha_i})$ has the form (a) as in Figure \ref{fig:4}. Since $e\geq 0$ then each crossing of $e$ has 
the form (b) as in Figure \ref{fig:4}. Hence, there exists an orientation of $L$ which makes it positive. Therefore, by [Ru2] $L$ is 
strongly quasipositive.
\end{proof}
\begin{figure}[H]
\centering
\def\svgwidth{0.50\columnwidth}
 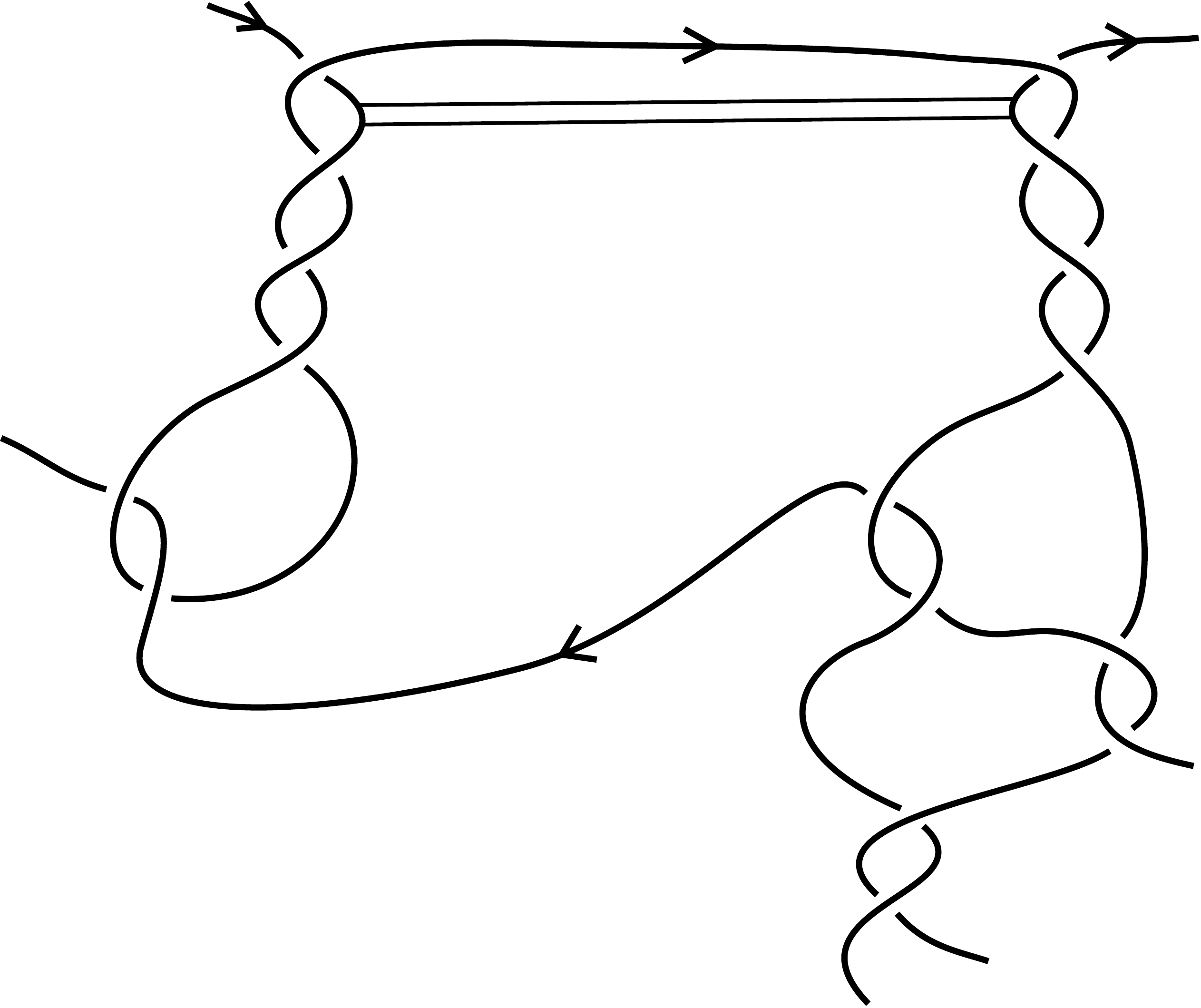
 \caption{Band move.}
\label{fig:5}
\end{figure}
\begin{proof}[Proof of Proposition \ref{Pthm3}] We prove this Proposition in the case where min$\{m_{i_0}, m_{i_0+1}\}=2$ because the general case is
done similarly.
With some orientation of $K$ , the $i_{0}^{th}$ and $({i_{0}+1})^{th}$ strands are oriented as in Figure \ref{fig:5}, which illustrates the case
$c_1^{i_0}>0$ and $c_1^{i_0+1}<0$.

Performing a band move as shown in Figure \ref{fig:5}, gives a punctured 2-sphere $A$ in $\mathbb{S}^3\times [0,1]$ such that 
$A\cap (\mathbb{S}^3\times \{1\})=K$ and $A\cap(\mathbb{S}^3\times \{0\})=K'\# L(\frac{q}{p})$ where 
$$K'=M(e; t_1,t_2,\cdots,t_{i_0-1}, t_{i_0+2}, \cdots, t_r)$$ and $L(\frac{q}{p})$ is the two-bridge link with rational slope
$$\frac{q}{p}=[2c_2^{i_0}+ 2c_2^{i_0+1}, 2c_3^{i_0+1},\cdots, 2c_{m_{i_0+1}}^{i_0+1}]$$ or 
$$\frac{q}{p}=[2c_2^{i_0+1}+ 2c_2^{i_0}, 2c_3^{i_0},\cdots, 2c_{m_{i_0}}^{i_0}].$$
Since $m_{i_0+1}$ and $m_{i_0}$ are even, $L(\frac{q}{p})$ bounds a Seifert surface $\Gamma$
of genus $$g(L(\frac{q}{p}))< max(\frac{m_{i_0}}{2}, \frac{m_{i_0+1}}{2}).$$
In $\mathbb{S}^3\times\{0\}$, $K'\# L(\frac{q}{p})$ bounds the boundary connected sum $S$ of $\Gamma$ with a minimal genus Seifert surface 
for $K'$. Therefore, taking the union of this surface $S$ with the punctured 2-sphere $A$ shows that
$$g_4(K)\leq g(K') + g(L(\frac{q}{p}))+1.$$
To complete the proof we compare the genus $g_4(K)$ and $g(K)$, and to compute the genus of $K$ and $K'$ we use [HM], we have three cases:

\textbf{First case:} If $e\neq 0$, then 

\begin{equation} \label{eq1}
\begin{split}
g(K) &=\frac{1}{2}(1+\sum_{i=1}^rm_i)\\
  &= \frac{1}{2}(1+\sum_{i\neq \{i_0, i_0+1\}}m_i) + \frac{m_{i_0}}{2} + \frac{m_{i_0+1}}{2}\\
  &= g(K') + \frac{m_{i_0}}{2} + \frac{m_{i_0+1}}{2}\\
  &= g(K') + max(\frac{m_{i_0}}{2}, \frac{m_{i_0+1}}{2})+1\\
  &> g(K')+g(L(\frac{q}{p}))+1\geq g_4(K).
\end{split}
\end{equation}
\textbf{Second case:} If $e= 0$ and ($c_1^1,c_1^2,\cdots, c_1^r)\neq \pm (1,-1,\cdots,1,-1)$ then 

\begin{equation} \label{eq1}
\begin{split}
g(K) &=\frac{1}{2}(-1+\sum_{i=1}^rm_i)\\
  &= \frac{1}{2}(-1+\sum_{i\neq \{i_0, i_0+1\}}m_i) + \frac{m_{i_0}}{2} + \frac{m_{i_0+1}}{2}\\
  &= g(K') + \frac{m_{i_0}}{2} + \frac{m_{i_0+1}}{2}\\
  &= g(K') + max(\frac{m_{i_0}}{2}, \frac{m_{i_0+1}}{2})+1\\
  &> g(K')+g(L(\frac{q}{p}))+1\geq g_4(K).
\end{split}
\end{equation}
\textbf{Third case:} If $e= 0$ and ($c_1^1,c_1^2,\cdots, c_1^r)= \pm (1,-1,\cdots,1,-1)$ then 
$g(K) =\frac{1}{2}(1+\sum_{i=1}^rm_i)-(p+1)$, where $p=min\{p_1,\cdots, p_r\}$, $p_i$ is the number of leading $2$'s in $t_i$ if $i$ is odd,
or the number of leading $-2$'s, if $i$ is even. Therefore, we have that
\begin{equation} \label{eq1}
\begin{split}
g(K) &=\frac{1}{2}(1+\sum_{i=1}^rm_i)-(p+1)\\
  &= \frac{1}{2}(1+\sum_{i\neq \{i_0, i_0+1\}}m_i)-(p+1) + \frac{m_{i_0}}{2} + \frac{m_{i_0+1}}{2}\\
  &\geq g(K') + \frac{m_{i_0}}{2} + \frac{m_{i_0+1}}{2}\\
  &\geq g(K') + max(\frac{m_{i_0}}{2}, \frac{m_{i_0+1}}{2})+1\\
  &> g(K')+g(L(\frac{q}{p}))+1\geq g_4(K).
\end{split}
\end{equation}
Since $g(K)>g_4(K)$ in all cases then $K$ is not strongly quasipositive.
\end{proof}

\end{document}